\documentclass[11pt]{article}
\usepackage{amsmath, amsthm, amssymb, amsfonts, latexsym,amscd}
\usepackage{authblk}
\usepackage{tikz}
\usepackage{mathtools}
\usepackage{xcolor}
\usetikzlibrary{shadings,patterns}
\newtheorem{theorem}{Theorem}[section]
\newtheorem{corollary}[theorem]{Corollary}
\newtheorem{lemma}[theorem]{Lemma}
\newtheorem{proposition}[theorem]{Proposition}
\newtheorem{definition}[theorem]{Definition}

\begin{document}

\title{ \bf On structural properties of some probable $R(3, 10)$-critical graphs}
\author[1]{Dinesh Pandey \thanks{dpandey@wlu.ca} \thanks {The first author’s research was funded by the Einwechter Centre for Supply Chain Management (a research centre funded by Mr. Dan Einwechter) at Wilfrid Laurier University.}}

\author[2]{Peruvemba Sundaram Ravi \thanks{pravi@wlu.ca}} 
\affil[1,  2] {Lazaridis School of Business and Economics,
             Wilfrid Laurier University, Waterloo, Ontario, N2L3C5, Canada.}
\date{}
\maketitle
\begin{abstract}
The Ramsey number $R(s, t)$ is the smallest positive integer $n$ such that every graph on $n$ vertices contains either a clique of size $s$ or an independent set of size $t$. An $R(s,t)$-critical graph is a graph on $R(s,t)-1$ vertices that contains neither a clique of size $s$ nor an independent set of size $t$. It is known that $40\leq R(3, 10)\leq 42$. We study the structure of a $R(3,10)$-critical graphs by assuming $R(3, 10)=42$. We show that if such a graph exists then its minimum degree and vertex connectivity are the same and is $6, 7$ or $8$. Then we find all the possible degree sequences of such graphs. Further, we show that if such a graph exists, then its diameter is either $2$ or $3$, and if it has diameter $2$ and minimum degree $6$, then it has only $21$ choices for its degree sequence.\\


\end{abstract}

\section{Introduction and Preliminaries}
The Ramsey number $R(s,t)$ is the smallest positive integer $n$ such that any graph on $n$ vertices either contains a clique of size $s$ or an independent set of size $t$. A graph on $R(s,t)-1$ vertices that contains neither a clique of size $s$ nor an independent set of size $t$ is called a $R(s,t)$-critical graph. Determining Ramsey numbers or establishing significant bounds for them poses a substantial challenge in combinatorics. It is easy to see that $R(1,t)=1$ for $t\geq 1$ and $R(2,t)=t$ for $t\geq 2$.  For $s, t \geq 3$, only a few Ramsey numbers have been established, namely $R(3, 3), R(3, 4), R(3, 5), (3, 6)$, $R(3, 7), R(3, 8)$, $R(3, 9), R(4, 4)$ and $R(4, 5)$. For some Ramsey numbers the upper and/or lower bounds are known. For comprehensive information on these bounds, we recommend the dynamic survey \cite{Radz}. Researchers are continuously working on refining these bounds.  $R(3, 10)$ is the smallest Ramsey number currently not known, although close bounds have been established.  

Our aim is to investigate the properties of $R(3,10)$-critical graphs and to use those properties to prove that $R(3,10)\leq 41$ without relying on computational tools. When we started work on this paper, it was known that $40 \leq R(3,10) \leq 42$. Therefore, in this paper, we initially assume that $R(3,10)=42$ and explore the structural properties of $R(3,10)$-critical graphs with $41$ vertices. These properties are presented in Section \ref{(3,10,41) graphs} of this paper. Subsequently, we became aware of work by Angelveit \cite{Angel} who used a computationally intensive approach to prove that $R(3,10) \leq 41$. We believe that the results presented in Section \ref{(3,10,41) graphs}  could be used in conjunction with Angeltveit's work, possibly to obtain a faster computationally intensive proof of the result that $R(3,10) \leq 41$. Furthermore, in Section  \ref{(3,10,40) graphs} - which was written after we became aware of Angelveit's work - we start with the assumption that $R(3,10)=41$ and employ our techniques (developed in Section \ref{(3,10,41) graphs}) to study potential $R(3,10)$-critical graphs with $40$ vertices, assuming their existence.  We believe that these techniques may also aid in the broader investigation of $R(3,t)$-critical graphs.

In this paper, we consider simple, finite, and undirected graphs. Let $G$ be a graph with vertex set $V(G)$ and edge set $E(G)$. For any vertex $v\in V(G)$, we denote the degree of $v$ in $G$ as $d_G(v)$ or simply $d(v)$. If $H$ is a subgraph of $G$ and $u\in V(H)$, then $d_H(u)$ indicates the degree of vertex $u$ in $H$. The notation $G-v$ refers to the graph obtained by removing vertex $v$ along with all edges incident to it. For two graphs $G$ and $H$, the union $G\cup H$ is the graph  having vertex set $V(G)\cup V(H)$ and the edge set $E(G)\cup E(H)$. The {\it neighbourhood} of vertex $v$, denoted $N(v)$, consists of all vertices adjacent to $v$. The {\it closed neighbourhood}, $N[v]=N(v)\cup \{v\}$. $G_v$ denotes the subgraph of $G$ induced by the vertices $V(G)\setminus N[v]$. We use $\delta(G)$ and $\Delta(G)$ to denote the minimum and maximum degree of $G$, respectively. The {\it independence number} of $G$ is denoted by $\alpha(G)$. A {\it vertex cut} in $G$ is a subset of $V(G)$ such that $G\setminus V(G)$ is disconnected. The {\it vertex connectivity} $\kappa(G)$ is the minimum number of vertices whose removal disconnects $G$ or reduces it to a single vertex. Similarly, the {\it edge connectivity} $\kappa'(G)$ is the minimum number of edges that must be removed to disconnect $G$. The {\it length} of a path in $G$ is the number of edges in the path. The {\it distance} $d(u,v)$ between two vertices $u$ and $v$ in $G$ is the length of the shortest path joining them. The {\it diameter} of $G$, denoted as $diam(G)$, is given by $diam(G)=\max\{d(u,v):u,v\in V(G)\}$. The {\it degree sequence} of $G$ is a list of its vertex degrees, typically arranged in descending or ascending order. In this paper, we present the degree sequences in descending order. The following fundamental result popularly known as the first theorem of Graph Theory, is frequently utilized in our analysis.

\begin{lemma}\label{first theorem}
In  any graph $G$, $\sum_{v\in V(G)}{d(v)}=2|E(G)|$.
\end{lemma}
 Let $G$ be a graph on $n$ vertices and $e$ edges which neither contains a clique of size $k$ nor an independent set of size $l$.  Then $G$ is said to be a $(k,l,n,e)$ graph or a $(k,l,n)$ graph or a $(k,l)$ graph. This notation is taken from \cite{Goed}. Also,
              $$e(k,l,n)=\min\{|E(G)|: G \mbox{ is a $(k,l,n)$ graph}\}.$$
 
\begin{lemma}[\cite{Goed}, Theorem3] \label {R(3,9) critical}
There is a unique $(3,9,35)$ graph and is $8$-regular.
\end{lemma}

We recall the Mantel's theorem, which gives an upper bound to the number of edges in a triangle free graph.
\begin{proposition} [{\bf Mantel's theorem}]\label{Mantel}
A triangle free graph on $n$ vertices contains at most $\frac{n^2}{4}$ edges.
\end{proposition}

\section{(3,10,41) graphs}\label{(3,10,41) graphs}
Assuming that $R(3,10)= 42$, we study the properties of $R(3,10)$-critical graphs with $41$ vertices without depending on computational tools. In this section $\Gamma$ is a $(3,10,41)$ graph.

\subsection{Minimum degree and vertex connectivity of $\Gamma$}
As $\Gamma$ is triangle free, no two neighbours of $v$ can be adjacent. Now since $\Gamma$ does not contain a $\bar{K_{10}}$, $\Delta(\Gamma)\leq 9$. As $|E(\Gamma)|\geq 172$ (see \cite{Goed}, Table 5), it follows from Lemma \ref{first theorem} that $\Gamma$ has a vertex of degree $9$. Thus $\Delta(\Gamma)=9$. In fact, it follows that $\Gamma$ has at least $16$ vertices of degree $9$. We establish a bound on the minimum degree of $\Gamma$.

\begin{lemma}\label {degree_sum_12}
 There is no vertex in $\Gamma$ such that the degree sum of two of its neighbours is less than or equal to $12$.
 \end{lemma}
\begin{proof}
Suppose $v\in V(\Gamma)$ and $v_1, v_2\in N(v)$ such that $d(v_1)+d(v_2)\leq 12$. Let $X=\{v_1, v_2\}$ and let $N[X]=N[v_1]\cup N[v_2]$. Then $|N[X]|\leq 13$, implying $|V(\Gamma) \setminus N[X]|\geq 28 = R(3,8)$. So $\Gamma \setminus N[X]$ contains an independent set $S$ of size 8. Note that no vertex of $\Gamma \setminus N[X]$ is adjacent to $v_1$ or $v_2$. So $S\cup \{v_1,v_2\}$ is an independent set of size 10 in $\Gamma$, which is a contradiction.
\end{proof}

\begin{theorem}
$\delta(\Gamma)\geq 6$.
\end{theorem}

\begin{proof}
Suppose $\Gamma$ has a vertex $v$ with $d(v)\leq 4$. Then $\Gamma_v$ is a graph on at least $36$ vertices. As $R(3,9)=36$, $\Gamma_v$ contains an independent set $S$ of size $9$. Thus $S$ together with $v$ forms an independent set of size 10 in $\Gamma$, which is a contradiction. This shows that $\delta(\Gamma)\geq 5$. Now it is sufficient to show that  $\Gamma$ has no vertex of degree $5$.

Suppose $\Gamma$ has a vertex $v$ of degree $5$. Then $\Gamma_v$ is a $(3,9,35)$ graph. So by Lemma \ref{R(3,9) critical}, $\Gamma_v$ is $8$-regular. Let $N(v)=\{v_1,v_2,v_3,v_4,v_5\}$. We show that non of $v_1, v_2,\ldots,v_5$ has degree more than $6$. Assume  $d(v_1)\geq 7$.  Note that neighbours of $v_1$ in $\Gamma_v$ are not adjacent to any of $v_2, v_3,v_4, v_5$ as their degrees now been exhausted. But $v_1$ has at least 6 neighbours in $\Gamma_v$ which is a contradiction as together with $v_2,v_3,v_4, v_5$ they form an independent set of size at least $10$ in $\Gamma$. This implies $d(v_i)\in \{5, 6\}$ for $i=1,2,3,4,5$. But this contradicts Lemma \ref {degree_sum_12}. Hence, $v$ can not be of degree $5$.
\end{proof}
\begin{lemma}\label{number_degree 6 vertices}
Number of vertices of degree $6$ in $\Gamma$ is at most $6$.
\end{lemma}

\begin{proof}
Let $v_1, v_2, \ldots, v_k$ be the $k$ vertices of degree $6$ in $\Gamma$. Then By Lemma \ref{degree_sum_12}, $N(v_i)$'s are disjoint. So $|V(\Gamma)|=41 \geq |N(v_1)\cup N(v_2)\cup\cdots \cup N(v_k)|=|N(v_1)|+|N(v_2)|+\cdots +|N(v_k)|=6k$. This implies $k\leq 6$.
\end{proof}

The following is a well known relation among vertex connectivity, edge connectivity and minimum degree of a simple graph.
\begin{proposition} \label{relation_Graph parameters}
For a simple graph $G$, $\kappa(G)\leq \kappa'(G)\leq \delta(G)$.
\end{proposition}

It is to observe that $\Gamma$ is a connected graph. A proof can be seen using contradiction in the following way: Suppose $H_1, H_2,\ldots, H_k$ are $k$ connected components of $\Gamma$ where $k\geq 2$.  Let $H_1=G_1$ and $\cup_{i=2}^k{H_i}=G_2$. Let $|V(G_1)|=n_1$. Then $|V(G_2)|=41-n_1$. For any value of $n_1\in\{1,2,\ldots, 40\}$, it can be easily checked that $\alpha(G_1)+\alpha(G_2)\geq 10$. This implies $\Gamma$ has an independent set of size at least $10$, which is a contradiction. In fact the connectivity of $\Gamma$ is much higher. We prove the following.


\begin{theorem} \label{vertex_connectivity}
 $\kappa(\Gamma)\geq 6$.
 \end{theorem}
 \begin{proof}
Suppose $\kappa(\Gamma)\leq 5$. Then there is a vertex cut of size $5$. Let $S=\{u_1,u_2,u_3,u_4,u_5\}$ be a vertex cut in $\Gamma$. Let $H_1, H_2,\ldots, H_k$ be $k$ connected components of $\Gamma \setminus S$ where $k\geq 2$.  Take $H_1=G_1$ and $\cup_{i=2}^k{H_i}=G_2$. Let $|V(G_1)|=n_1$. Then $|V(G_2)|=n_2=36-n_1$. Without loss of generality assume that $n_1\leq n_2$. Then $1\leq n_1\leq 18.$

\noindent{\bf Case 1:} $n_1=1$\\ 
Suppose $V(G_1)=\{v\}$.  Then $d_{\Gamma}(v)\leq 5$, a contradiction.\\
{\bf Case 2:} $n_1=2$\\
Suppose $V(G_1)=\{v_1,v_2\}$. In this case $n_2=34$ and hence $G_2$ contains an independent set $S_2$ of size $8$. This implies $v_1\sim v_2$ otherwise $S_2$ together with $v_1$ and $v_2$ form an independent set of size $10$ in $\Gamma$. As $\Gamma$ is triangle free, non of $u_1,u_2,\ldots, u_5$ is adjacent to both $v_1$ and $v_2$. This implies $d(v_1)+d(v_2)\leq 7$ and hence either $d(v_1)<4$ or $d(v_2)<4$, a contradiction.\\
{\bf Case 3:} $3 \leq n_1 \leq 8$.\\
In this case, $n_2\geq 28$ and hence $G_2$ contains an independent set $S_2$ of size $8$. Also since $n_1\geq 3$, $G_1$ contains an independent set $S_1$ of size $2$. $S_1\cup S_2$ is an independent set of size $10$ in $\Gamma$, a contradiction.\\
{\bf Case 4:} $9\leq n_1\leq 18$\\
In this case $G_1$ contains an independent set of size $4$ and $G_2$ contains an independent set of size $6$ which collected together gives an independent set of size $10$ in $\Gamma$, a contradiction.\\
This completes the proof.
 \end{proof}
 
 \begin{theorem}
$\kappa(\Gamma)=\delta(\Gamma)$.
 \end{theorem}
 
 \begin{proof}
 As $|V(\Gamma)|=41$, $\Gamma$ is not $9$-regular. So $6\leq \delta(\Gamma)\leq 8$.\\
 
 \noindent{\bf Case I:} $\delta(\Gamma)=6$\\
  By Proposition \ref{relation_Graph parameters} and Theorem \ref{vertex_connectivity}, it follows that $\kappa(\Gamma)=6$.\\
 
\noindent {\bf Case II:}\label{min_deg_8} $\delta(\Gamma)=8$\\
 Suppose $\Gamma$ has a vertex cut of size $7$. Let $S=\{u_1,u_2,\ldots, u_7\}$ be a vertex cut in $\Gamma$.  Let $H_1, H_2,\ldots, H_k$ be $k$ connected components of $\Gamma \setminus S$ where $k\geq 2$.  Take $H_1=G_1$ and $\cup_{i=2}^k{H_i}=G_2$. Let $|V(G_1)|=n_1$. Then $|V(G_2)|=n_2=34-n_1$. Without loss of generality assume that $n_1\leq n_2$. Then $1\leq n_1\leq 17.$ Now we can consider the following six cases: $(i) n_1=1, (ii) n_1=2 (iii) 3\leq n_1\leq 6 (iv) 7\leq n_1\leq 11 (v) 12\leq n_1 \leq 16$ and $(vi) n_1=18$. In each of these cases we get a contradiction by giving similar arguments as in the proof of Theorem \ref{vertex_connectivity}. So $\Gamma$ can not have a vertex cut of size $7$. This implies $\kappa(\Gamma)\geq 8$ and hence  $\kappa(\Gamma)=8$.\\
 
 \noindent{\bf Case III:}$\delta(\Gamma)=7$\\
 Similar arguments as in Case II proves that $\kappa(\Gamma)=7$.
 \end{proof}
 
\subsection{Possible degree sequences of $\Gamma$}
In this section we find all the possible degree sequences of $\Gamma$ provided it exists. We saw that $\Delta(\Gamma)=9$ and $\delta(\Gamma)\geq 6$. So $6\leq d(v)\leq 9$ $\forall v\in V(\Gamma)$. We represent the degree sequence of $\Gamma$ by the $4$-tuple $(a, b, c, d)$ where $a, b, c$ and $d$ are the vertices of degrees $9, 8, 7$ and $6$, respectively. As $|V(\Gamma)|=41$, all the vertices of $\Gamma$ can not be of odd degrees. So $\Gamma$ must have a vertex of degree $6$ or $8$. So $6\leq \delta(\Gamma)\leq 8$. We divide the graph $(3,10,41)$ graphs into two classes. One when $\delta(\Gamma)>6$ and the other when $\delta(\Gamma)=6$. 

\subsubsection{Degree sequences of $\Gamma$ when $\delta(\Gamma)>6$}

It is known that $172\leq |E(\Gamma)|\leq 184$.  First suppose $\Gamma$ is a $(3, 10, 41, 172)$-graph.  As $\delta(\Gamma)> 6$, $d=0$. So we have 

\begin{equation*}
9a+8b+7c=344, 
\end{equation*}
\begin{equation*}
a+b+c=41,
\end{equation*}
and
\begin{equation*}
a=\alpha. 
\end{equation*}
From Lemma \ref{first theorem}, it follows that $16\leq \alpha \leq 28$. Also for $16\leq \alpha \leq 28$,

\[ rank
\left(
\begin{array}{cccc}
 9 & 8  & 7 & 344 \\
1  & 1  & 1 & 41  \\
 1 & 0  & 0 & \alpha   
\end{array}
\right)=
rank \left(
\begin{array}{cccc}
 9 & 8  & 7  \\
1  & 1  & 1  \\
 1 & 0  & 0    
\end{array}
\right)= 3.
\]
So the above system of equation has a unique solution for each $\alpha$. The solutions (in the form $(a, b, c, d)$) are $(16, 25, 0, 0), (17, 23, 1, 0 ), (18, 21, 2, 0), (19, 19, 3, 0), (20, 17, 4, 0), (21, 15, 5, 0)$,\\ $(22, 13, 6, 0), (23, 11, 7, 0), (24, 9, 8, 0), (25, 7, 9, 0), (26, 5, 10, 0), (27, 3, 11, 0), (28, 1, 12, 0)$. Each of these solutions is a possibility for the degree sequence of $\Gamma$ when it is a $(3, 10, 41, 172)$- graph. Similarly we can obtain the possible degree sequence of every $(3, 10, 41, e)$- graph for $173\leq e\leq 184$. The possible degree sequences of a $(3, 10, 41, e)$-graph for $172\leq e\leq 184$ are listed in Table \ref{no degree 6 vertex}.

\begin{table}[h!]
\begin{center}
\begin{tabular}{|c| l |}
\hline
$e$& \hskip 5 cm possible degree sequences\\
\hline
172& $(16, 25, 0, 0), (17, 23, 1, 0 ), (18, 21, 2, 0), (19, 19, 3, 0), (20, 17, 4, 0), (21, 15, 5, 0), (22, 13, 6, 0), $\\ 
&$ (23, 11, 7, 0), (24, 9, 8, 0), (25, 7, 9, 0), (26, 5, 10, 0), (27, 3, 11, 0), (28, 1, 12, 0)$\\
\hline
173& $(18, 23, 0, 0), (19, 21, 1, 0), (20, 19, 2, 0), (21, 17, 3, 0), (22, 15, 4, 0), (23, 13, 5, 0), (24, 11, 6, 0),$\\
 &$ (25, 9, 7, 0), (26, 7, 8, 0), (27, 5, 9, 0), (28, 3, 8, 0), (29, 1, 9, 0)$\\
\hline
174 & $(20, 21, 0, 0), (21, 19, 1, 0), (22, 17, 2, 0), (23, 15, 3, 0), (24, 13, 4, 0), (25, 11, 5, 0), (26, 9, 6, 0),$\\
 & $ (27, 7, 7, 0), (28, 5, 8, 0), (29, 3, 9, 0), (30,1, 10, 0)$\\
\hline 
175 & $(22, 19, 0, 0), (23, 17, 1, 0), (24, 15, 2, 0), (25, 13, 3, 0), (26, 11, 4, 0), (27, 9, 5, 0), (28, 7, 6, 0),$\\
 &$ (29, 5, 7, 0), (30, 3, 8, 0), (31,1, 9, 0)$\\
\hline
176 & $(24, 17, 0, 0), (25, 15, 1, 0), (26, 13, 2, 0), (27, 11, 3, 0), (28, 9, 4, 0), (29, 7, 5, 0), (30, 5, 6, 0),$ \\
&$ (31, 3, 7, 0), (32, 1, 8, 0)$\\
\hline 
177 & $(26, 15, 0, 0), (27, 13, 1, 0), (28, 11, 2, 0), (29, 9, 3, 0), (30, 7, 4, 0), (31, 5, 5, 0), (32, 3, 6, 0)$,\\
 & $ (33, 1, 7, 0)$ \\
\hline 
178 & $(28, 13, 0, 0), (29, 11, 1, 0), (30, 9, 2, 0), (31, 7, 3, 0), (32, 5, 4, 0), (33, 3, 5, 0), (34, 1, 6, 0)$\\
\hline 
179 & $(30, 11, 0, 0), (31, 9, 1, 0), (32, 7, 2, 0), (33, 5, 3, 0), (34, 3, 4, 0), (35, 1, 5, 0)$ \\
\hline 
180 & $(32, 9, 0, 0), (33, 7, 1, 0), (34, 5, 2, 0), (35, 3, 3, 0), (36, 1, 4, 0)$\\
\hline
181 & $(34, 7, 0, 0), (35, 5, 1, 0), (36, 3, 2, 0), (37, 1, 3, 0)$\\
\hline 
182 & $(36, 5, 0, 0), (37, 3, 1, 0), (38, 1, 2, 0)$ \\
\hline 
183 & $(38, 3, 0, 0), (39, 1, 1, 0)$ \\
\hline
184 & $(40, 1, 0, 0)$\\
\hline
\end{tabular}
\end{center}
\label{default}
\caption{possible degree sequences of a $(3, 10, 41, e)$-graph of minimum degree at least 7} \label{no degree 6 vertex}
\end{table}

\subsubsection{degree sequences of $\Gamma$ when $\delta(\Gamma)=6$}
Since $\delta(\Gamma)=6$, by Lemma \ref{number_degree 6 vertices} we have $1\leq d\leq 6$. For each value of $d$ the degree sequences can be obtained as above. First suppose $d=1$. \\
If $\Gamma$ is a $(3, 10, 41, 172)$-graph, then we have 
\begin{equation*}
9a+8b+7c=338,
\end{equation*}
\begin{equation*}
a+b+c=40,
\end{equation*}
and
\begin{equation*}
a=\alpha. 
\end{equation*}
By Lemma \ref{first theorem}, it follows that $18\leq \alpha \leq 29$.  This system of equation has a unique solution for each $\alpha$ and every solution $(a, b, c)$ represent the degree sequence $(a, b, c, 1)$ of $\Gamma$. These degree sequences are listed in Table \ref{1 degree 6 vertex}. Similarly the degree sequences for $d=2, 3, 4, 5$ and $6$ can be obtained which are listed in Tables \ref{1 degree 6 vertex}, \ref{2 degree 6 vertex}, \ref{3 degree 6 vertex}, \ref{4 degree 6 vertex}, \ref{5 degree 6 vertex} and \ref{6 degree 6 vertex}, respectively.

\begin{table}[h!]
\begin{center}
\begin{tabular}{|c| l |}
\hline
$e$& \hskip 5 cm possible degree sequences\\
\hline
172& (18, 22, 0, 1), (19, 20, 1, 1), (20, 18, 2, 1), (21, 16, 3, 1), (22, 14, 4, 1), (23, 12, 5, 1), \\
      & (24, 10, 6, 1), (25, 8, 7, 1), (26, 6, 8, 1), (27, 4, 9, 1), (28, 2, 10, 1), (29, 0, 11, 1) \\ 
\hline
173& (20, 20, 0, 1), (21, 18, 1, 1), (22, 16, 2, 1), (23, 14, 3, 1), (24, 12, 4, 1), (25, 10, 5, 1) \\
      & (26, 8, 6, 1), (27, 6, 7, 1), (28, 4, 8, 1), (29, 2, 9, 1), (30, 0, 10, 1)\\
\hline
174 & (22, 18, 0, 1), (23, 16, 1, 1), (24, 14, 2, 1), (25, 12, 3, 1), (26, 10, 4, 1) (27, 8, 5, 1)\\
       & (28, 6, 6, 1), (29, 4, 7, 1), (30, 2, 8, 1), (31, 0, 9, 1)\\
\hline 
175 & (24, 16, 0, 1), (25, 14, 1, 1), (26, 12, 2, 1), (27, 10, 3, 1), (28, 8, 4, 1), (29, 6, 5, 1),\\ 
        & (30, 4, 6, 1), (31, 2, 7, 1), (32, 0, 8, 1)\\
\hline
176 & (26, 14, 0, 1), (27, 12, 1, 1), (28, 10, 2, 1), (29, 8, 3, 1), (30, 6, 4, 1), (31, 4, 5, 1),\\
       &(32, 2, 6, 1), (33, 0, 7, 1) \\
\hline 
177 & (28, 12, 0, 1), (29, 10, 1, 1), (30, 8, 2, 1), (31, 6, 3, 1), (32, 4, 4, 1), (33, 2, 5, 1),\\
       &(34, 0, 6, 1)\\
\hline 
178 & (30, 10, 1, 1), (31, 8, 8, 1), (32, 6, 2, 1), (33, 4, 3, 1), (34, 2, 4, 1), (35, 0, 5, 1)\\
\hline 
179 &  (32, 8, 0, 1), (33, 6, 1, 1), (34, 4, 2, 1), (35, 2, 3, 1), (36, 0, 4, 1)\\
\hline 
180 & (34, 6, 0, 1), (35, 4, 1, 1), (36, 2, 2, 1), (37, 0, 3, 1) \\
\hline
181 & (36, 4, 0, 1), (37, 2, 1, 1), (38, 0, 2, 1) \\
\hline 
182 &  (38, 2, 0, 1), (39, 0, 1, 1)\\
\hline 
183 & (40, 0, 0, 1) \\
\hline

\end{tabular}
\end{center}
\label{default}
\caption{\small possible degree sequences of a $(3, 10, 41, e)$-graph containing exactly one vertex of degree $6$} \label{1 degree 6 vertex}
\end{table}

\begin{table}[h!]
\begin{center}
\begin{tabular}{|c| l |}
\hline
$e$& possible degree sequences\\
\hline
172& (20, 19, 0, 2), (21, 17, 1, 2), (22, 15, 2, 2), (23, 13, 3, 2), (24, 11, 4, 2), (25, 9, 5, 2),\\
       &(26, 7, 6, 2), (27, 5, 7, 2), (28, 3, 8, 2), (29, 1, 9, 2)  \\ 
\hline
173& (22, 17, 0, 2), (23, 15, 1, 2), (24, 13, 2, 2), (25, 11, 3, 2), (26, 9, 4, 2), (27, 7, 5, 2),\\
      & (28, 5, 6, 2), (29, 3, 7, 2), (30, 1, 8, 2)\\
\hline
174 & (24, 15 , 0 , 2), (25, 13, 1, 2), (26, 11, 2, 2), (27, 9, 3, 2), (28, 7, 4, 2), (29, 5, 5, 2)\\
       & (30, 3, 6, 2), (31, 1, 7, 2)\\
\hline 
175 & (26, 13, 0, 2), (27, 11, 1, 2), (28, 9, 2, 2), (29, 7, 3, 2), (30, 5, 4, 2), (31, 3, 5, 2)\\
       & (32, 1, 6, 2)\\
\hline
176 & (28, 11, 0, 2), (29, 9, 1, 2), (30, 7, 2, 2), (31, 5, 3, 2), (32, 3, 4, 2), (33, 1, 5, 2)\\
\hline 
177 & (30, 9, 0, 2), (31, 7, 1, 2), (32, 5, 2, 2), (33, 3, 3, 2), (34, 1, 4, 2)\\
\hline 
178 & (32, 7, 0, 2), (33, 5, 1, 2), (34, 3, 2, 2), (35, 1, 3, 2)\\
\hline 
179 & (34, 5, 0, 2), (35, 3, 1, 2), (36, 1, 2, 2)\\
\hline 
180 &(36, 3, 0, 2), (37, 1, 1, 2)\\
\hline
181 & (38, 1, 0, 2)\\
\hline 
\end{tabular}
\end{center}
\label{default}
\caption{possible degree sequences of a $(3, 10, 41, e)$-graph containing two vertices of degree $6$}  \label{2 degree 6 vertex}
\end{table}

\begin{table}[h!]
\begin{center}
\begin{tabular}{|c| l |}
\hline
$e$& possible degree sequences\\
\hline
172& (22, 16, 0, 3), (23, 14, 1, 3), (24, 12, 2, 3), (25, 10, 3, 3), (26, 8, 4, 3), (27, 6, 5, 3),\\
      & (28, 4, 6, 3), (29, 2, 7, 3), (30, 0, 8, 3)\\ 
\hline
173& (24, 14, 0, 3), (25, 12, 1, 3), (26, 10, 2, 3), (27, 8, 3, 3), (28, 6, 4, 3), (29, 4, 5, 3)\\
      & (30, 2, 6, 3) (31, 0, 7, 3)\\
\hline
174 & (26, 12, 0, 3), (27, 10, 1, 3), (28, 8, 2, 3), (29, 6, 3, 3), (30, 4, 4, 3), (31, 2, 5, 3) \\
       & (32, 0, 6, 3)\\
\hline 
175 & (28, 10, 0, 3), (29, 8, 1, 3), (30, 6, 2, 3), (31, 4, 3 ,3), (32, 2, 4, 3), (33, 0, 5, 3) \\
\hline
176 & (30, 8, 0, 3), (31, 6, 1, 3), (32, 4, 2, 3), (33, 2, 3, 3), (34, 0, 4, 3)\\
\hline 
177 & (32, 6, 0, 3), (33, 4, 1, 3), (34, 2, 2, 3), (35, 0, 3, 3)\\
\hline 
178 & (34, 4, 0, 3), (35, 2, 1, 3), (36, 0, 2, 3)\\
\hline 
179 & (36, 2, 0, 3), (37, 0, 1, 3)\\
\hline 
180 & (38, 0, 0, 3)\\
\hline
\end{tabular}
\end{center}
\label{default}
\caption{possible degree sequences of a $(3, 10, 41, e)$-graph containing three vertices of degree $6$}  \label{3 degree 6 vertex}
\end{table}

\begin{table}[h!]
\begin{center}
\begin{tabular}{|c| l |}
\hline
$e$& possible degree sequences\\
\hline
172& (24, 13, 0, 4), (25, 11, 1, 4), (26, 9, 2, 4), (27, 7, 3, 4), (28, 5, 4, 4), (29, 3, 5, 4) \\
      & (30, 1, 6, 4)\\ 
\hline
173& (26, 11, 0, 4), (27, 9, 1, 4), (28, 7, 2, 4), (29, 5, 3, 4), (30, 3, 4, 4), (31, 1, 5, 4)\\
\hline
174 & (28, 9, 0, 4), (29, 7, 1, 4), (30, 5, 2, 4), (31, 3, 3, 4), (32, 1, 4, 4)\\
\hline 
175 & (30, 7, 0, 4), (31, 5, 1, 4), (32, 3, 2, 4), (33, 1, 3, 4)\\
\hline
176 & (32, 5, 0, 4), (33, 3, 1, 4), (34, 1, 2, 4)\\
\hline 
177 & (34, 3, 0, 4), (35, 1, 1, 4)\\
\hline 
178 & (36, 1, 0, 4)\\
\hline 

\end{tabular}
\end{center}
\label{default}
\caption{possible degree sequences of a $(3, 10, 41, e)$-graph containing four vertices of degree $6$}  \label{4 degree 6 vertex}
\end{table}

\begin{table}[h!]
\begin{center}
\begin{tabular}{|c| l |}
\hline
$e$& possible degree sequences\\
\hline
172& (26, 10, 0, 5), (27, 8, 1, 5), (28, 6, 2, 5), (29, 4, 3, 5), (30, 2, 4, 5), (31, 0, 5, 5)\\ 
\hline
173& (28, 8, 0, 5), (29, 6, 1, 5), (30, 4, 2, 5), (31, 2, 3, 5), (32, 0, 4, 4)\\
\hline
174 & (30, 6, 0, 5), (31, 4, 1, 5), (32, 2, 2, 5), (33, 0, 3, 5)\\
\hline 
175 & (32, 4, 0, 5), (33, 2, 1, 5), (34, 0, 2, 5)\\
\hline
176 & (34, 2, 0, 5), (35, 0, 1, 5)\\
\hline 
177 & (36, 0, 0, 5)\\
\hline 
\end{tabular}
\end{center}
\label{default}
\caption{possible degree sequences of a $(3, 10, 41, e)$-graph containing five vertices of degree $6$}  \label{5 degree 6 vertex}
\end{table}

\begin{table}[h!]
\begin{center}
\begin{tabular}{|c| l |}
\hline
$e$& possible degree sequences\\
\hline
172& (28, 7, 0, 6), (29, 5, 1, 6), (30, 3, 2, 6), (31, 1, 3, 6)\\ 
\hline
173& (30, 5, 0, 6), (31, 3, 1, 6), (32, 1, 2, 6)\\
\hline
174 & (32, 3, 0, 6), (33, 1, 1, 6)\\
\hline 
175 & (34, 1, 0, 6)\\
\hline
\end{tabular}
\end{center}
\label{default}
\caption{possible degree sequences of a $(3, 10, 41, e)$-graph containing six vertices of degree $6$}  \label{6 degree 6 vertex}
\end{table}

\subsection{diameter of $\Gamma$}
In this section we show that  $diam(\Gamma)$ is either $2$ or $3$. We analyse the structure of $\Gamma$ when it is of diameter $2$ and minimum degree $6$. When $diam(\Gamma)=2$ and $\delta(\Gamma)=6$, it is shown that the number of possible degree sequences of $\Gamma$ is only $21$.

\begin{theorem}
$2\leq diam(\Gamma)\leq 3$.
\end{theorem}

\begin{proof}
The left inequality is obvious as $\Gamma$ is not a complete graph. To show the right inequality, it is sufficient to show that the distance between any two vertices of $\Gamma$ is at most $3$. Suppose there exists $u, v\in V(\Gamma)$ such that $d(u,v)\geq 4$. Then $N(u)\cap N(v)=\emptyset $  and no vertex in $N(u)$ is adjacent to a vertex in $N(v)$, otherwise the distance between $u$ and $v$ will be $3$ or less. As $\delta(\Gamma)\geq 6$, $N(\{u,v\})$ forms an independent set of size at least $12$, which is a contradiction. 
\end{proof}

  \begin{figure}[h!]
 \begin{center}
 \begin{tikzpicture}
 \filldraw (0,0) node[above]{$v$} circle [radius=.5mm] (1,-1)node [right]{$v_4$}  circle [radius=.5mm] (2,-1)node [right]{$v_5$}  circle [radius=.5mm] (3,-1)node [right]{$v_6$}  circle [radius=.5mm] (-1,-1)node [left]{$v_3$}  circle [radius=.5mm] (-2,-1)node [left]{$v_2$}  circle [radius=.5mm] (-3,-1)node [left]{$v_1$} circle [radius=.5mm];
 \draw (-1,-1)--(0,0)--(1,-1) (-2,-1)--(0,0)--(2,-1) (-3,-1)--(0,0)--(3,-1);
 \draw (0, -3)node {$G_v$} circle [x radius=3, y radius=1];
 \end{tikzpicture}
 \caption {Structure of the graph $G$}
 \end{center}
 \end{figure}
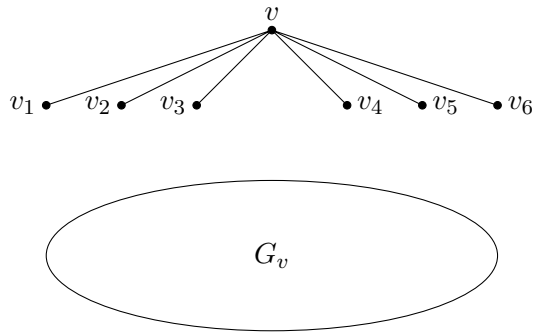
 
 \begin{theorem}\label {degree 2}
If $diam(\Gamma)=2$, then $\Gamma$ has at most two vertices of degree $6$. Furthermore, if there are two vertices of degree $6$, then they must be adjacent.
\end{theorem}

\begin{proof}
Let $v$ be a vertex of degree $6$ in $\Gamma$. As $diam(\Gamma)=2$, it follows from Lemma \ref{degree_sum_12} that there is no vertex of degree $6$ in $\Gamma_v$. From Lemma \ref{degree_sum_12}, it also follows that $N(v)$ contains at most one vertex of degree $6$. Hence, $\Gamma$ can have at most two vertices of degree $6$ and if there exist two, then they must be adjacent.
\end{proof}

\begin{lemma}\label{$H2v_8$}
Let $v$ be a vertex of degree $6$ in $\Gamma$. Then 
$$20\leq |\{u\in V(\Gamma_v): d_{\Gamma_v}(u)=8\}|\leq 24.$$
\end{lemma}

\begin{proof}
Note that $\Gamma_v$ induces a $(3,9,34)$-graph. The left inequality holds from the fact that $e(3,9,34)=129$ (see \cite{Goed}, Table 4). For convenience let us name the set $\{u\in \Gamma_v: d_{\Gamma_v}(u)=8\}$ as $(\Gamma_v,8)$.To prove the right inequality, suppose $|(\Gamma_v,8)|\geq 25.$ As $|N(v)|=6$, so by Pigeonhole principal, there is a vertex $w$ in $N(v)$ which has at least $5$ neighbours in $(\Gamma_v,8)$. These $5$ neighbours of $w$ together with the remaining $5$ vertices in $N(v)$ forms an independent set of size $10$ in $\Gamma$, which is a contradiction.
\end{proof}

\begin{definition}
Let $G$ be a graph and $v\in V(G)$. For $0\leq i\leq \Delta(G)$, we define
 $$H_{2i}(v)=\{u\in \Gamma_v: |N(u)\cap N(v)|=i\}.$$
\end{definition}

\begin{proposition}
Let $diam(\Gamma)=2$ and $v$ be a vertex of degree $6$ in $\Gamma$. Then $H_{20}(v)=\emptyset= H_{2i}(v)$ for $i\geq 4.$
\end{proposition}
\begin{proof}
$H_{20}(v)=\emptyset$ because $diam(\Gamma)=2$. Now observe that ${\{u\in \Gamma_v: d_{\Gamma_v}(u)\leq 5\}=\emptyset}$, otherwise $\Gamma_v\setminus N[u]$ contains an independent set $S$ of size $8$, which together with $\{u,v\}$ forms an independent set of size $10$ in $\Gamma.$ Let $i\geq 4$ and let $w\in H_{2i}(v)$. Then $d_{\Gamma_v}(w)\geq 6$. Hence $d_{\Gamma}(w)= 6+i \geq 10$, a contradiction. Hence $H_{2i}(v)=\emptyset$ for $ i\geq 4.$
\end{proof}
 
\begin{proposition}
Let $diam(\Gamma)=2$ and $v$ be a vertex of degree $6$ in $\Gamma$. Then $20\leq |H_{21}(v)|\leq 24$.
\end{proposition}

\begin{proof}
As $diam(\Gamma)=2$, $(\Gamma_v,8)\subseteq H_{21}(v)$ and hence the left inequality follows from Lemma \ref{$H2v_8$}.  If $|H_{21}(v)|\geq 25$, then by Pigeon hole principle, there is a vertex $w$ in $N(v)$ which has at least $5$ neighbours in $H_{21}(v)$. These five neighbours of $w$ together with $N(v)\setminus\{w\}$ contains an independent set of size $10$ in $\Gamma$, which is a contradiction. Therefore $|H_{21}(v)|\leq 24$.
\end{proof}

Note that here $|H_{21}(v)|+|H_{22}(v)|+|H_{23}(v)|= 34$ and there are at most $48$ edges between $N(v)$ and $\Gamma_v$.  This implies if $|H_{21}(v)|=20$, then $|H_{23}(v)|=0$ and so $|H_{22}|=14$. Similarly we can obtain all the possible values of the triplet $(|H_{21}(v)|, |H_{22}(v)|, |H_{23}(v)|)$ . These values are listed in Table \ref{edge distribution} along with the number of edges $|E(N(v),\Gamma_v)|$ between $N(v)$ and $\Gamma_v$ in each case.

\begin{table}[h!]
\begin{center}
\begin{tabular}{|c|c|c|c|}
\hline
$|H_{21}(v)|$&$|H_{22}(v)|$&$|H_{23}(v)|$&$|E(N(v),\Gamma_v)|$\\
\hline
20&14&0&48\\
\hline
21&13&0&47\\
\hline
21&12&1&48\\
\hline
22&12&0&46\\
\hline
22&11&1&47\\
\hline
22&10&2&48\\
\hline
23&11&0&45\\
\hline
23&10&1&46\\
\hline
23&9&2&47\\
\hline 
23&8&3&48\\
\hline
24&10&0&44\\
\hline
24&9&1&45\\
\hline
 24&8&2&46\\
\hline
24&7&3&47\\
\hline 
24&6&4&48\\
\hline
\end{tabular}
\end{center}
\label{default}
\caption{All possible values of $(|H_{21}(v)|, |H_{22}(v)|, |H_{23}(v)|)$ } \label{edge distribution}
\end{table}

In Table \ref{edge distribution}, we got that $44\leq |E(N(v), \Gamma_v)|\leq 48$. The next Lemma describes the contribution to the degree sequence of $\Gamma$ by the vertices in $N[v]$.
\begin{lemma}\label{contri_N[v]}
Let $diam(\Gamma)=2$ and $v$ be a vertex of degree $6$ in $\Gamma$. Let $(a_1, b_1, c_1, d_1)$ be the contribution to the degree sequence of $\Gamma$ by the vertices in $N[v]$. Then the following hold. 
\begin{enumerate}
\item \label{edge_44} if $|E(N(v), \Gamma_v)|= 44$ then $(a_1, b_1, c_1, d_1)\in \{(4, 1, 0, 2), (4, 0, 2, 1), (3, 2, 1, 1), (2, 4, 0, 1)\}$,
\item if $|E(N(v), \Gamma_v)|= 45$ then $(a_1, b_1, c_1, d_1)\in \{(5, 0, 0, 2), (4, 1, 1, 1), (3, 3, 0, 1)\}$,
\item if $|E(N(v), \Gamma_v)|= 46$ then $(a_1, b_1, c_1, d_1)\in \{(5, 0, 1, 1), (4, 2, 0, 1)\}$,
\item if $|E(N(v), \Gamma_v)|= 47$ then $(a_1, b_1, c_1, d_1)=(5, 1, 0, 1)$,
\item  if $|E(N(v), \Gamma_v)|= 48$ then $(a_1, b_1, c_1, d_1)= (6,0, 0, 1)$.
\end{enumerate}
\end{lemma}

\begin{proof}
Suppose $|E(N(v), \Gamma_v)|= 44$. Let $N(v)=\{v_1,v_2,\ldots, v_6\}$. Then we have 
\begin{equation*}
9a_1+8b_1+7c_1+6d_1=56,
\end{equation*}
and
\begin{equation*}
a_1+b_1+c_1+d_1=7.
\end{equation*}

From Lemma \ref{first theorem}, it follows that $2\leq a_1\leq 4$. Also from Theorem \ref{degree 2}, we have $1\leq d_1\leq 2$. With these restrictions on $a_1$ and $d_1$, it can be easily checked that the only solutions of the above system in terms of $(a_1, b_1, c_1, d_1)$ are $(4, 1, 0, 2), (4, 0, 2, 1), (3, 2, 1, 1)$ and $(2, 4, 0, 1)$. This completes the proof of {\it \ref{edge_44}}. Similar arguments can be given to prove the other parts.
\end{proof}
 We now count the vertices of degrees $6, 7, 8$ and $9$ in $\Gamma_v$. Clearly $\sum_{u\in \Gamma_v}{d(u)}\leq 34\times 9=306$ and $\sum_{u\in \Gamma_v}{d(u)}\geq \sum_{u\in \Gamma_v}{d_{\Gamma_v}(u)}+44\geq 129\times 2 +44= 302$. So we have $302\leq \sum_{u\in \Gamma_v}{d(u)}\leq 306$.

\begin{lemma}\label{contri_Gamma_v}
Let $diam(\Gamma)=2$ and $v$ be a vertex of degree $6$ in $\Gamma$. Let $(a, b, c, d)$ be the number of vertices of degrees $ 9, 8, 7$ and $6$ respectively in $\Gamma_v$. Then the following hold.
\begin{enumerate}
\item \label{sum_302} if $\sum_{u\in \Gamma_v}{d(u)}=302$, then $(a, b, c, d)\in \{(32,1,0,1), (32,0,2,0), (31,2,1,0),(30,4,0,0)\}$,
\item if $\sum_{u\in \Gamma_v}{d(u)}=303$, then $(a, b, c, d)\in \{(33,0,0,1), (32,1,1,0), (31,3,0,0)\}$,
\item if $\sum_{u\in \Gamma_v}{d(u)}=304$, then $(a, b, c, d)\in (33,0,1,0), (32,2,0,0)\}$,
\item if $\sum_{u\in \Gamma_v}{d(u)}=305$, then $(a, b, c, d)=(33,1,0,0)$,
\item If $\sum_{u\in \Gamma_v}{d(u)}=306$, then $(a, b, c, d)=(34,0,0,0)$.
\end{enumerate}
\end{lemma}

\begin{proof}
  Suppose $\sum_{u\in \Gamma_v}{d(u)}=302$. As $a, b, c$ and $d$ be vertices in $\Gamma_v$ having degrees $9, 8, 7$ and $6$ in $\Gamma$. We have

 \begin{equation*}
9a+8b+7c+6d=302,
\end{equation*} 
\begin{equation*}
a+b+c+d=34.
\end{equation*}
By Lemma \ref{first theorem} it follows that  $30\leq a\leq 32$. Now it can be easily checked that, 
\begin{itemize}
\item $a=30 \implies (b, c, d)=(4, 0, 0)$
\item $a=31 \implies (b, c, d)=(2, 1, 0)$ \mbox{and} 
\item $a= 30 \implies (b, c, d)\in \{(1, 0, 1), (0, 2, 0)\}$.
\end{itemize} 
This proves {\it \ref{sum_302}}. Similar arguments can be given to prove the other parts.
\end{proof}

\begin{theorem}\label {degree sequence_6}
Let $diam(\Gamma)=2$ and $\Gamma$ has a vertex of degree $6$. Then the degree sequence of $\Gamma$ is one of the following $21$ sequences. \\
$(40, 0,0,1), (39, 0, 1, 1), (38, 2,0,1),  (38, 0, 2, 1), (38, 1, 0, 2), (37, 2, 1, 1), (37, 1, 1, 2), (37, 0, 3, 1)$,\\ $(36, 4, 0, 1), (36, 3, 0, 2), (36, 2, 2, 1), (36, 1, 2, 2), (36, 0, 4, 1), (35, 4, 1, 1), (35, 3, 1, 2), (35, 2, 3, 1)$,\\ $(34, 6, 0, 1), (34, 5, 0, 2),   (34, 4, 2, 1), (33, 6, 1, 1), (32, 8, 0, 1)$. 
\end{theorem}

\begin{proof}

We write the degree sequence of $\Gamma$ as $(a, b, c, d)$ where $a, b, c$ and $d$ are the  number of vertices of degree $9, 8, 7$ and $6$, respectively.  As $N[v]$ and $\Gamma_v$ partitions $V(\Gamma)$, from Lemma \ref{contri_N[v]} and Lemma \ref{contri_Gamma_v} we have $(a, b, c, d)=(a_1, b_1, c_1, d_1)+(a_2, b_2, c_2, d_2)$ where $a_1, b_1, c_1, d_1$ are the number of vertices of degrees $9, 8, 7, 6$ respectively in $N[v]$ and $a_2, b_2, c_2, d_2$ are the number of vertices of degrees $9, 8, 7, 6$ respectively in $\Gamma_v$.  We have $|E(N(v), \Gamma_v)|\in \{44, 45,46, 47, 48\}$. As $\Gamma_v$ is a $(3, 9, 34)$-graph, $|E(\Gamma_v)|\geq 129$. If $|E(\Gamma_v)|\geq 132$, then $\sum_{u\in \Gamma_v}{d(u)}\geq 264+44=308$, a contradiction. So $|E(\Gamma_v)|\in \{129, 130, 131\}$. If
\begin{itemize}
\item $|E(N(v), \Gamma_v)|=44$, then by Lemma \ref{contri_N[v]}, $(a_1, b_1, c_1, d_1)\in \{(4, 1, 0, 2), (4, 0, 2, 1), (3, 2, 1, 1)$, \\ $(2, 4, 0, 1)\}$:= $S_{N[v]}$. Also in this case  $\sum_{u\in \Gamma_v}{d(u)}\in \{302, 304, 306\}$. So by Lemma \ref{contri_Gamma_v}, $(a_2, b_2, c_2, d_2)\in \{(32,1,0,1), (32,0,2,0), (31,2,1,0),(30,4,0,0), (33,0,1,0), (32,2,0,0),\\
(34,0,0,0)\}$:=$S_{\Gamma_v}$. So we get the degree sequences of $\Gamma$ 
\begin{align*}
(a, b, c, d)&\in\{(a_1, b_1, c_1, d_1)+ (a_2, b_2, c_2, d_2): (a_1, b_1, c_1, d_1)\in S_{N[v]}, (a_2, b_2, c_2, d_2)\in S_{\Gamma_v},\\
                &\;\;\;\;\;\; d_1+d_2\leq 2\}\\
 		 &=\{(36, 1, 2, 2), (35, 3, 1, 2), (34, 5, 0, 2), (37, 1, 1, 2), (36, 3, 0, 2), (38, 1, 0, 2),\\
		 &\;\;\;\;\;\; (36, 0, 4, 1), (35, 2, 3, 1), (34, 4, 2, 1), (37, 0, 3, 1), (36, 2, 2, 1), (38, 0, 2, 1),\\
		 &\;\;\;\;\;\; (35, 4, 1, 1), (37, 2, 1, 1), (33, 6, 1, 1), (32, 8, 0, 1), (34, 6, 0, 1), (36, 4, 0, 1)\}.
\end{align*}
By giving similar arguments the following can be shown.
\item If $|E(N(v), \Gamma_v)|=45$, then the set of possible degree sequences of $\Gamma$ is
\begin{align*}
\{& (37, 1, 1, 2), (36, 3, 0, 2), (38, 1, 0, 2), (36, 2, 2, 1), (35, 4, 1, 1), (37, 2, 1, 1),\\
 & (34, 6, 0, 1), (36, 4, 0, 1)\}
\end{align*}
\item $|E(N(v), \Gamma_v)|=46$, then the set of possible degree sequences of $\Gamma$ is\\
$\{(38, 0, 2, 1), (37, 2, 1, 1), (39, 0, 1, 1), (36, 4, 0, 1), (38, 2, 0, 1)\}$
\item $|E(N(v), \Gamma_v)|=47$, then the set of possible degree sequences of $\Gamma$ is $\{(38, 2, 0, 1)\}$.
\item $|E(N(v), \Gamma_v)|=48$, then the set of possible degree sequences of $\Gamma$ is $\{(40, 0,0,1)\}$
\end{itemize}
Now the result follows by taking the union of all the above possible degree sequences of $\Gamma$.
\end{proof}

\section{(3,10,40) graphs} \label{(3,10,40) graphs}
In this section, we assume that $R(3,10)=41$. Then there exists an $R(3,10)$-critical graph on $40$ vertices. Let $\Omega$ be a $(3,10,40)$ graph. We draw some conclusions about $\Omega$  (assuming that it exists). As $e(3,10,40) \geq 161$, it follows that
$\Delta(\Omega)=9$. Also $\delta(\Omega)\geq 40-R(3,9)=4$.

\begin{lemma}\label {degree_sum_11}
 There is no vertex in $\Omega$ such that the degree sum of two of its neighbours is less than or equal to $11$.
 \end{lemma}
The proof follows by similar arguments as in the proof of Theorem \ref{degree_sum_12}. Let $\Omega_v$ be the graph $\Omega\setminus N[v]$.

\begin{lemma}
Let $v\in V(\Omega)$ and $v_1,v_2\in N(v)$. Then $|N(v_1)\cup N(v_2)|\geq11$. 
\end{lemma}
\begin{theorem}
$diam(\Omega)=2$ or $3$.
\end{theorem}
\begin{proof}
Clearly $diam(\Omega)\geq 2$. It is sufficient to prove that $d(u,v)\neq 4$ for any $u, v\in V(\Omega)$. Suppose $u,v\in V(\Omega)$ such that $d(u,v)=4$. Then $N(u)\cup N(v)$ forms an independent set in $\Omega$.  If any of $u$ or $v$ has degree $6$ or more then $|N(u)\cup N(v)|\geq 10$. Hence both $u$ and $v$ has degree at most $5$. So $\Omega\setminus  (N[u]\cup N[v])$ is a graph on at least $29$ vertices and hence contains an independent set of $S$ of size $8$ (because $R(3,8)=28$). $S\cup \{u,v\}$ forms an independent set of size $10$ in $\Omega$, which is a contradiction. 
\end{proof}

\begin{corollary}
Let $diam(\Omega)=2$ and $v$ be a vertex of degree $6$ in $\Omega$. Then there can be at most $2$ vertices of degree $6$ in $N(v)$. 
\end{corollary}

\begin{proof}
Let $N(v)=\{v_1,v_2,v_3,v_4,v_5,v_6\}$. We first show that if two vertices in $N(v)$ has degree $6$, then their neighbours in $\Gamma_v$ are disjoint. Suppose $d(v_1)=d(v_2)=6$ and $(N(v_1)\cap N(v_2))\cap V(\Gamma_v)\neq \emptyset$, then $|N[v_1]\cap N[v_2]|\leq 12$. So $G\setminus (N[v_1]\cup N[v_2])$ is a graph on at least $28$ vertices and hence contains an independent set $S$ of size $8$. $S\cup\{v_1,v_2\}$ is an independent set of size $10$ in $\Gamma$, a contradiction. Hence $(N(v_1)\cap N(v_2))\cap V(\Gamma_v)= \emptyset$.\\

Suppose $N(v)$ contains three vertices of degree $6$. and let $d(v_1)=d(v_2)=d(v_3)=6$. Then $N(v_1)\cup N(v_2)\cup N(v_3)$ has $15$ vertices in $\Gamma_v$. Let $H= V(\Gamma_v)\setminus (N(v_1)\cup N(v_2)\cup N(v_3))$. If any of $v_4, v_5$ or $v_6$ has $7$ neighbours in $H$, then they together with $\{v_1,v_2,v_3\}$ form an independent set of size $10$ in $\Gamma$, a contradiction. So suppose none of $v_4, v_5$ and $v_6$ has $7$ neighbours in $H$. As $diam(\Gamma)=2$ and $|H|=18$, each of $v_4, v_5$ and $v_6$ has $6$ neighbours in $H$. Also no neighbour of $v_4$ is adjacent to either $v_5$ or $v_6$, otherwise $v_5$ or $v_6$ gets $7$ neighbours in $H$. So $(N(v_4)\cap H)\cup \{v_1,v_2,v_3,v_5,v_6\}$ is an independent set of size $11$ in $\Gamma$, a contradiction.
\end{proof}

If $\Omega$ is $9$-regular, then its diameter must be $2$. So the results holding for diameter $2$ (3, 10, 40)- graphs also holds for the probable $9$-regular (3, 10, 40) graphs. Let $\Omega_i(v)=\{x\in V(\Omega): d(x,v)=i\}$. With this notation, $\Omega_0(v)=\{v\}, \Omega_1(v)=N(v)$ and for $i\geq 4$, $\Omega_i(v)=\emptyset$. Next lemma discusses the possible sizes of $\Omega_i(v)$ for $i=2$ and $3$ when $d(v)=4$.

\begin{lemma}\label {set_sizes}
Let $v\in V(\Omega)$ with $d(v)=4$. Then $19\leq |\Omega_2(v)|\leq 24$ and $11\leq |\Omega_3(v)|\leq 17$.
\end{lemma}

\begin{proof}
If $|\Omega_3(v)|\geq 18$, then $\Omega_3(v)$ contains an independent set $S$ of size $6$ which together with $N(v)$ forms an independent set of size $10$, a contradiction. So $|\Omega_3(v)|\leq 17$. Suppose $|\Omega_3(v)|\leq 10$. Then $|\Omega_2(v)|\geq 25$. Let $H$ be the subgraph induced by $\Omega_2(v)\cup \Omega_3(v)$ in $\Omega$. Then $H$ is a  $(3,9,35)$ graph and hence is $8$-regular from Lemma \ref{R(3,9) critical}. Hence every vertex of $\Omega_2(v)$ is adjacent to exactly one vertex in $N(v)$. Let $N(v)=\{v_1, v_2, v_3, v_4\}$. Since $|\Omega_2(v)|\geq 25$, there exists a vertex (say) $v_1$ in $N(v)$ with $d(v_1)\geq 8$ i.e. $N(v_1)\cap \Omega_2(v)\geq 7$. $N(v_1)\cap \Omega_2(v)$ together with $\{v_2,v_3,v_4\}$ forms an independent set of size $10$ in $\Omega$, a contradiction. Hence $|\Omega_3(v)|\geq 11$ and we get $10\leq |\Omega_3(v)|\leq 17$.

    As $|\Omega_2(v)|+|\Omega_3(v)|=35$, it follows that $18\leq |\Omega_2(v)|\leq 24$. Suppose $|\Omega_2(v)|=18$. Then $d(v_1)+d(v_2)+d(v_3)+d(v_4)=22$. By Lemma \ref{degree_sum_11}, $d(v_1)+d(v_2)\geq 11$ and $d(v_3)+d(v_4)\geq 11$. ThIs implies $d(v_1)+d(v_2)= 11$ and $d(v_3)+d(v_4)= 11$. Let $d(v_1)\leq d(v_2)$ and $d(v_3)\leq d(v_4)$.  Then $d(v_1)\leq 5$ and $d(v_3)\leq 5$, which contradicts Lemma \ref {degree_sum_11}. Hence $|\Omega_2(v)|\geq 19$ and so $19\leq |\Omega_2(v)|\leq 24$.
 \end{proof}

Denote $\Omega_v^i:=\{x\in V(\Omega_v): |N(x) \cap N(v)|=i\}$.

\begin{theorem} 
If $diam(\Omega)=2$, then $\delta(\Omega)\geq 6$.
\end{theorem}
\begin{proof}
It is known that $\delta(\Omega)\geq 4$. Suppose there exists a vertex $v$ of degree $4$ in $\Omega$. As $diam(\Omega)=2$, $|\Omega_2(v)|=35$. This contradicts Lemma \ref{set_sizes}.\\

Now suppose there exists a vertex $v$ of degree $5$ in $\Omega$. Then $\Omega_v=G\setminus N[v]$ is a $(3,9,34)$ graph. So $\Omega_v^1\geq 20$. If $\Omega_v^1> 20$, then there exists a vertex $u$ in $N(v)$ such that $N(u)\cap N(\Omega_v^1)\geq 6$. So $N(u)\cap N(\Omega_v^1)$ together with the remaining $4$ vertices of $N(v)$ forms an independent set of size $10$, a contradiction. So $|\Omega_v^1|=20$. This implies $|E(N(v)\rightarrow \Omega_v)|\geq 20+28=48$, a contradiction. So $\Omega$ can not have a vertex of degree $5$. This implies $\delta(\Omega)\geq 6$.
\end{proof}

\begin{lemma} \label{min_connectivity_Omega}
If $diam(\Omega)=2$, then $\kappa(\Omega)\geq 6$.
\end{lemma}
\begin{proof}
Suppose $\kappa(\Gamma)\leq 5$. Then there is a vertex cut of size $5$. Let $S=\{u_1,u_2,u_3,u_4,u_5\}$ be a vertex cut in $\Gamma$. Let $H_1, H_2,\ldots, H_k$ be $k$ connected components of $\Gamma \setminus S$ where $k\geq 2$.  Take $H_1=G_1$ and $\cup_{i=2}^k{H_i}=G_2$. Let $|V(G_1)|=n_1$. Then $|V(G_2)|=n_2=35-n_1$. Without loss of generality assume that $n_1\leq n_2$. Then $1\leq n_1\leq 17.$ 
Rest of the proof is similar to the proof of Lemma \ref{vertex_connectivity}, by considering the five cases, $n_1=1, n_1=2, 3\leq n_1\leq 7, 8\leq n_1\leq 12$ and $13\leq n_1\leq 17$.
\end{proof}

\begin{theorem}
If $diam(\Omega)=2$, then $\kappa(\Omega)=\delta(\Omega)$.
\end{theorem}
\begin{proof}
$6\leq \delta(\Omega)\leq 9$. Consider the following cases.\\ 
 \noindent{\bf Case I:} $\delta(\Omega)=6$\\
  By Proposition \ref{relation_Graph parameters} and Lemma \ref{min_connectivity_Omega}, it follows that $\kappa(\Omega)=6$.\\
  
 \noindent {\bf Case II:}\label{min_deg_8} $\delta(\Omega)=9$\\
 Suppose $\Omega$ has a vertex cut of size $8$. Let $S=\{u_1,u_2,\ldots, u_8\}$ be a vertex cut in $\Omega$.  Let $H_1, H_2,\ldots, H_k$ be $k$ connected components of $\Omega \setminus S$ where $k\geq 2$.  Take $H_1=G_1$ and $\cup_{i=2}^k{H_i}=G_2$. Let $|V(G_1)|=n_1$. Then $|V(G_2)|=n_2=32-n_1$. Without loss of generality assume that $n_1\leq n_2$. Then $1\leq n_1\leq 16.$ Now consider the following seven cases: $(i) n_1=1, (ii) n_1=2 (iii) 3\leq n_1\leq 4 (iv) n_1=5 (v) 6\leq n_1\leq 8 (vi)9 \leq n_1 \leq 13$ and $(vii) 14\leq n_1=16$. \\
 
\noindent If $n_1=1$, then $G_1$ is a a single vertex graph $v$ and $d_{\Omega}(v)=8$, a contradiction.\\

\noindent If $n_1=2$, then suppose $V(G_1)=\{v_1,v_2\}$. In this case $n_2=34$ and hence $G_2$ contains an independent set $S_2$ of size $8$. This implies $v_1\sim v_2$ otherwise $S_2$ together with $v_1$ and $v_2$ form an independent set of size $10$ in $\Omega$. As $\Omega$ is triangle free, non of $u_1,u_2,\ldots, u_8$ is adjacent to both $v_1$ and $v_2$. So $d(v_1)+d(v_2)\leq 10$ and hence either $d(v_1)<5$ or $d(v_2)<5$, a contradiction.\\

\noindent If $3\leq n_1\leq 4$, then $G_1$ contains an independent set of size $2$ and $G_2$ contains an independent set of size 8 and hence $\Omega$ contains an independent set of size 10, a contradiction.\\

\noindent If $n_1=5$, then by Proposition \ref {Mantel}, number of edges in $G_1$ is at most 6. Let $V(G_1)=\{v_1, v_2, v_3, v_4, v_5\}$. Each vertex of $S$ is adjacent to at most two vertices in $V(G)$, otherwise $G_1$ is not $R(3,3)$-critical. So $d_{\Omega}(v_1)+d_{\Omega}(v_2)+d_{\Omega}(v_3)+d_{\Omega}(v_4)+d_{\Omega}(v_5)\leq 28$. This implies $d_{\Omega}(v_i)\leq 5$ for some $1\leq i\leq 5$, a contradiction.\\

\noindent For rest of the cases it can easily be shown (as in the case $3\leq n_1\leq 4$) that $\Omega$ contains an independent set of size 10 to get a contradiction. Thus $\kappa(\Omega)\geq 9$ and hence by Proposition \ref{relation_Graph parameters}, $\kappa(\Omega)=9=\delta(\Omega)$.\\

\noindent {\bf Case III:} $\delta(\Omega)=8$\\
 Similar arguments as in Case II can be given by taking the cut set $\{u_1, u_2,\ldots, u_7\}$ and considering the cases $n_1=1, n_1=2, 3\leq 5, 6\leq n_1\leq 10, 11 \leq n_1\leq 13, 14\leq n_1\leq 16$ to conclude that $\kappa(\Omega)=8$.
 
 \noindent {\bf Case IV:} $\delta(\Omega)=7$\\
 Similar arguments as in Case II can be given by taking the cut set $\{u_1, u_2,\ldots, u_6\}$ and considering the cases $n_1=1, n_1=2, 3\leq 6, 7\leq n_1\leq 11, 12 \leq n_1\leq 16$ to conclude that $\kappa(\Omega)=7$.
 \end{proof}
 
 \begin{theorem}
 Let $diam(\Omega)=2$ and $S$ be a smallest cut set in $\Omega$. Then $S=N(v)$ for some $v\in V(\Omega)$. 
 \end{theorem}
 
 \begin{proof}
Let $H_1, H_2,\ldots, H_k$ be $k$ connected components of $\Omega \setminus S$ where $k\geq 2$.  W.L.O.G. assume that $|V(H_1)|\leq |V(H_i)|$ for $2\leq i\leq k$. Take $H_1=G_1$ and $\cup_{i=2}^k{H_i}=G_2$. Let $|V(G_1)|=n_1$. As $S\leq 9$, $|V(G_2)|=n_2=31-n_1$. So $1\leq n_1\leq 15.$ To prove the result, we need to show that $n_1=1$.\\
Suppose $n_1\geq 2$. As $|S|=\kappa(\Omega)\leq \delta(\Omega)$, $6\leq |S|\leq 9$. First suppose $|S|=9$. Consider the cases $2 \leq n_1\leq 5, 6\leq n_1\leq 8, 9\leq n_1\leq 13$ and  $14\leq n_1\leq 15$. If $2\leq n_1\leq 5$, then  by Proposition \ref {Mantel}, number of edges in $G_1$ is at most 6 and no vertex in $S$ is adjacent to more than two vertices in $G_1$. So $\sum_{v_i\in V(G_1)}{d_{\Gamma}(v_i)}\leq 30$. This implies $G_1$ contains a vertex of degree (in $\Omega$) less than $9$, which is a contradiction. If $6\leq n_1\leq 8$, then $G_1$ contains an independent set of size $3$ and $G_2$ contains an independent set of size $7$, a contradiction. Similarly contradiction arises when $9\leq n_1\leq 13$ and  $14\leq n_1\leq 15$. Analogous  arguments can be given when $|S|=8,7$ or $6$ to get a contradiction. This completes the proof.
 \end{proof}
 
 \section{Some nearly equivalent problems}
We know the degree sequences of all probable $(3,10,41)$ graphs. There can be more than one graph with the same degree sequence. Can a method be developed to get an upper bound on the number of non-isomorphic graphs with the same degree sequence? For example in \cite {Goed}, it is shown that the $8$-regular graph on $35$ vertices is unique. If it is known that there are $r$ non-isomorphic graphs with the degree sequence $D_n$, can a method be developed to construct all of them? These problems could be related to the problem of finding a Ramsey number $R(s,t)$.

\vskip 1cm
\noindent {\bf Authors' Note:} When we began working on this problem, it was known that $40\leq R(3,10)\leq 42$. We sought to prove that $R(3,10)\leq 41$ without the use of computational tools. We subsequently became aware of Angelveit's paper (arXiv December 2023; Electronic Journal of Combinatorics November 2025). This paper provided a computationally intensive proof that $R(3,10)\leq 41$. In light of this result, we added Section \ref{(3,10,40) graphs} to our paper and have temporarily suspended the project while we consider possible ways to extend and reorient this work.
\end{document}